\documentclass[12pt]{article}
\input amssymb.sty
\input amssym.def
\input amssym
\input epsf

\newtheorem{theorem}{Theorem}[section]
\newtheorem{lemma}[theorem]{Lemma}
\newtheorem{proposition}[theorem]{Proposition}
\newtheorem{corollary}[theorem]{Corollary}
\newtheorem{conjecture}[theorem]{Conjecture}
\newtheorem{definition}[theorem]{Definition}

\begin{document}
\newcommand{\be}{\begin{equation}}
\newcommand{\ee}{\end{equation}}
\newcommand{\bt}{\begin{theorem}}
\newcommand{\et}{\end{theorem}}
\newcommand{\bd}{\begin{definition}}
\newcommand{\ed}{\end{definition}}
\newcommand{\bp}{\begin{proposition}}
\newcommand{\ep}{\end{proposition}}
\newcommand{\bl}{\begin{lemma}}
\newcommand{\el}{\end{lemma}}
\newcommand{\bc}{\begin{corollary}}
\newcommand{\ec}{\end{corollary}}
\newcommand{\bcon}{\begin{conjecture}}
\newcommand{\econ}{\end{conjecture}}
\newcommand{\la}{\label}
\newcommand{\Z}{{\Bbb Z}}
\newcommand{\R}{{\Bbb R}}
\newcommand{\Q}{{\mathbb Q}}
\newcommand{\C}{{\Bbb C}}
\newcommand{\hra}{\hookrightarrow}
\newcommand{\lra}{\longrightarrow}
\newcommand{\lms}{\longmapsto}

\begin{titlepage}
\title{Hidden Hodge symmetries and Hodge correlators}
\author{A.B. Goncharov }
\date{\it To Don Zagier for his 60-th birthday}

\end{titlepage}
\maketitle


\section{Hidden Hodge symmetries}


There is a well known parallel between 
Hodge and \'etale theories, still incomplete and
  rather mysterious:

\begin{center}\begin{tabular}{| c | c | c |}
\bf{$l$-adic \'Etale  Theory}&{\bf Hodge Theory}\\
\hline\hline  
Category of $l$-adic & 
Abelian category ${\cal M}{\cal H}_\R$ \\
Galois modules&of real mixed Hodge strucrures\\
\hline
Galois group & Hodge Galois group  $G_{\rm Hod}:= $\\
${\rm Gal}(\overline \Q/\Q)$ 
&  Galois group of the category ${\cal M}{\cal H}_\R$ \\
\hline
${\rm Gal}(\overline \Q/\Q)$ acts on 
$H^*_{\rm et}(\overline X, \Q_l)$, & $H^*(X(\C), \R)$ 
has a functorial \\
where $X$ is a variety over $\Q$& real mixed Hodge structure\\
\hline
\'etale site & ??\\
\hline
${\rm Gal}(\overline \Q/\Q)$ acts on the \'etale site, and thus  &??\\
on categories of \'etale sheaves on $X$, e.g. &??\\
on the category of $l$-adic perverse sheaves& ??\\
\hline
${\rm Gal}(\overline \Q/\Q)$-equivariant perverse sheaves&Saito's Hodge sheaves\\
\hline
\end{tabular}\end{center}
The current absense of the ``Hodge site'' was emphasized by 
A.A. Beilinson  \cite{B}.

\paragraph{The Hodge Galois group.} A weight $n$ {\it pure real Hodge structure} 
is a real vector space $H$ together with a decreasing filtration $F^\bullet H_\C$ 
on its complexification satisfying
$$
H_\C = \oplus_{p+q=n}F^p H_\C\cap \overline {F^q} H_\C.
$$
A real Hodge structure is a direct sum  of pure ones. 
The category ${\cal P}{\cal H}_\R$ real Hodge structures  
 is equivalent to 
the category of representations of the real algebraic group 
$\C_{\C/\R}^*$. The group of complex points 
of $\C_{\C/\R}^*$ is $\C^*\times \C^*$; the complex conjugation interchanges 
the factors. 

A real {\it mixed Hodge structure} is given by a real vector space $H$ equipped with 
the weight  filtration $W_\bullet H$ and the Hodge filtration 
$F^\bullet H_\C$ 
of its complexification, such that the Hodge filtration induces on 
${\rm gr}^W_nH$ a  weight $n$ 
real Hodge structure. 
The category ${\cal M}{\cal H}_\R$ 
of real mixed Hodge structures 
is an abelian rigid tensor category. There is a fiber functor to the category of 
real vector spaces
$$
\omega_{\rm Hod}: 
{\cal M}{\cal H}_\R \lra {\rm Vect}_\R, \qquad {H} \lra \oplus_n{\rm gr}^W_{n}{H}.
$$
The  Hodge Galois group is a real algebraic group given by automorphisms 
of the fiber functor:
$$
{G}_{\rm H}: = {\rm Aut}^\otimes\omega_{\rm Hod}.
$$ 
The fiber functor provides a canonical equivalence of  categories
$$
\omega_{\rm Hod}: {\cal M}{\cal H}_\R \stackrel{\sim}{\lra} G_{\rm Hod}-\mbox{modules}.
$$ 
The Hodge Galois group is a semidirect product of the unipotent radical $U_{\rm Hod}$ 
and  $\C_{\C/\R}^*$:
\be \la{HODD}
0 \lra U_{\rm Hod} \lra {G}_{\rm Hod} \lra {\C^*}_{\C/\R} \lra 0, 
\qquad {\C^*}_{\C/\R}\hra {G}_{\rm Hod}.
\ee
The projection ${G}_{\rm Hod} {\to} {\C^*}_{\C/\R}$ is provided by the 
inclusion of the category of real Hodge structurs to the category of  mixed 
real Hodge structures. 
The splitting $s: {\Bbb G}_m \to {G}_{\rm Hod}$ is provided by the functor $\omega_{\rm Hod}$. 

The complexified Lie algebra of $U_{\rm Hod}$ 
has {\it canonical} generators $G_{p,q}$, $p,q\geq 1$, 
satisfying the only relation $\overline G_{p,q} = - G_{q,p}$, 
defined in \cite{G1}. For the subcategory of Hodge-Tate structures 
they were defined in \cite{L}. Unlike similar but different 
Deligne's generators \cite{D}, they behave nicely in families. 
So to define an action of the group $ {G}_{\rm Hod}$ 
one needs to have an  action of the subgroup 
$\C^*_{\C/\R}$ and, in addition to this, an  action  of a single operator 
$$
G:= \sum_{p,q\geq 1} G_{p,q}. 
$$

\paragraph{The twistor Galois group.}
Denote by $\C^*$ the real algebraic group with the group of complex points $\C^*$. 
The extension induced from (\ref{HODD})  
by the diagonal embedding $\C^* \subset \C_{\C/\R}^*$ is the {\it twistor Galois group}.  
It is a semidirect product of the groups $U_{\rm Hod}$ and $\C^*$. 
\be \la{HODD2}
0 \lra U_{\rm Hod} \lra {G}_{\rm T} \stackrel{\longleftarrow}{\lra}  {\C^*} \lra 0. 
\ee
It is not difficult to prove
\bl
The category of representations of ${G}_{\rm T}$ is equivalent to the category of 
mixed twistor structures defined by Simpson \cite{Si2}. 
\el

We suggest the following fills the $??$-marks in the 
dictionary related the Hodge and Galois. Below 
$X$ is a smooth projective complex algebraic variety.

\bcon \la{1.1.11.2} 
There exists a functorial homotopy 
action of the twistor Galois group $G_{\rm T}$ by $A_\infty$-equivalences of 
an $A_\infty$-enhancement of the derived category of perverse sheaves on 
$X$ such that 
the category of equivariant objects is 
equivalent to Saito's category  real mixed Hodge sheaves.\footnote{We want to have a 
natural construction of the action first, and get 
Saito's category  real mixed Hodge sheaves as a consequence, not the other way around.} 
\econ

Denote by $D^b_{\rm sm}(X)$ the category of 
smooth complexes of sheaves on $X$, i.e. 
complexes of sheaves on $X$ whose cohomology are local systems. 

\bt \la{6.21.11.1}
There exists a functorial for pull-backs 
 homotopy 
action of the twistor Galois group $G_{\rm T}$ by $A_\infty$-equivalences of 
an $A_\infty$-enhancement of the category $D^b_{\rm sm}(X)$.
\et 

The action of the subgroup $\C^*$ is not algebraic. It arises from Simpson's action of 
 $\C^*$ on semisimple local systems \cite{Si1}.  
The action of the Lie algebra of the unipotent radical $U_{\rm Hod}$ 
is determined by a collection of numbers, which we call 
the {\it Hodge correlators for semisimple 
local systems}. 
Our construction uses 
the theory of harmonic bundles \cite{Si1}. The Hodge correlators  can be interpreted as 
 correlators for a certain Feynman integral. 
This Feynman integral is probably responsible for the ``Hodge site''. 

For the trivial 
local system the construction was carried out  in  \cite{G2}. 
A more general construction for curves, 
involving the constant sheaves and delta-functions, 
 was carried out  in  \cite{G1}.

In the case when $X$ is the universal modular curve, the Hodge correlators 
contain the special values $L(f,n)$ of 
weight $k\geq 2$ modular forms for $GL_2(\Q)$ outside of the critical strip 
-- it turns out that 
the simplest Hodge correlators in this case coincide with the 
Rankin-Selberg integrals 
for the non-critical special values 
$L(f,k+n)$, $n \geq 0$ -- the case $k=2, n=0$ is discussed 
in detail in \cite{G1}.

\section{Hodge correlators for local systems} \la{3.2Har}
\subsection{An action of $G_T$ on the `` minimal model''of ${\cal D}_{\rm sm}(X)$.} 
Tensor products of irreducible local systems 
are semisimple local systems. 
The {\it category of harmonic bundles } ${\rm Har}_X$ 
is   the graded category whose objects are 
semi-simple local systems on $X$ and their shifts, and morphisms are given by 
graded vector spaces
\be \la{Harhom}
{\rm Hom}^\bullet_{{\rm Har}_X}(V_1, V_2):= H^\bullet(X, V^\vee_1\otimes V_2).
\ee
 
Here is our main result.
\bt \la{1.5.11.1}
There is a homotopy action of the twistor Galois group $G_T$ by 
$A_\infty$-equivalences of the graded category ${\rm Har}_X$, such that 
the action of the subgroup $\C^*$ is given by Simpson's action of $\C^*$ 
on semi-simple local systems. 
\et

This immediately implies Theorem \ref{6.21.11.1}. Indeed, 
given a small $A_\infty$-category ${\cal A}$, there is 
a functorial constraction of the 
triangulated envelope ${\rm Tr}({\cal A})$ of ${\cal A}$, 
the smallest triangulated category containing ${\cal A}$. 
Since ${\cal D}^b_{\rm sm}(X)$ is generated as a triangulated category 
by semi-simple local systems, the category ${\rm Tr}({\rm Har}_X)$ 
is equivalent to ${\cal D}^b_{\rm sm}(X)$ as a triangulated category, and thus is 
an $A_\infty$-enhancement of the latter. 
On the other hand, the action of the group $G_T$ from 
Theorem \ref{1.5.11.1} extends by functoriality 
to the action on ${\rm Tr}({\rm Har}_X)$. 

Below we recall what are  $A_\infty$-equivalences of DG categories 
and then define the corresponding data in our case. 

\subsection{$A_\infty$-equivalences of DG categories}
\paragraph{The Hochshild cohomology of a small  dg-category ${\cal A}$.} 
Let ${\cal A}$ be a small dg category. 
Consider 
 a bicomplex whose $n$-th column is 
\be \la{HCCH}
\prod_{[X_i]}{\rm Hom}\Bigl({\cal A}(X_0, X_1)[1] \otimes {\cal A}(X_1, X_2)[1] 
\otimes \ldots \otimes {\cal A}(X_{n-1}, X_n)[1], {\cal A}(X_{0}, X_n)[1]\Bigr),
\ee
where the product is over isomorphism classes $[X_i]$ of objects of the category ${\cal A}$. 
The vertical differential $d_1$ in the bicomplex is given by the differential on the 
tensor product of complexes. The horisontal one $d_2$ is the degree $1$ map 
provided by the composition 
$$
{\cal A}(X_{i}, X_{i+1})\otimes {\cal A}(X_{i+1}, X_{i+2}) \lra {\cal A}(X_{i}, X_{i+2}).
$$
Let ${\rm HC}^*({\cal A})$
be the total complex of this bicomplex. 
Its cohomology are the Hochshild   cohomology ${\rm HH}^*({\cal A})$ of ${\cal A}$. 
Let ${\rm Fun}_{A_{\infty}}({\cal A}, {\cal A})$ be the space of $A_{\infty}$-functors 
from ${\cal A}$ to itself.  Lemma \ref{6.27.11.2} 
can serve as a definition of $A_{\infty}$-functors considered modulo homotopy equivalence.

\bl \la{6.27.11.2} One has 
\be \la{6.27.11.3}
H^0{\rm Fun}_{A_{\infty}}({\cal A}, {\cal A})= {\rm HH}^0({\cal A}).
\ee
\el
Indeed, a cocycle in ${\rm HC}^0({\cal A})$ is the same thing as an $A_{\infty}$-functor. 
Coboundaries corresponds to the homotopic to zero functors.

\paragraph{The cyclic homology of a small rigid dg-category ${\cal A}$.}  
Let  $(\alpha_0 \otimes ... \otimes \alpha_m)_{\cal C}$ be 
the projection of 
$\alpha_0 \otimes ... \otimes \alpha_m$ to the coinvariants of the cyclic shift.
So, if $\overline \alpha:= {\rm deg}\alpha$, 
$$
(\alpha_{0} \otimes ... \otimes \alpha_{m})_{\cal C} = (-1)^{\overline \alpha_{m} (\overline \alpha_{0}
+ ... + \overline \alpha_{m-1})}(\alpha_{1} \otimes ... \otimes \alpha_{m}\otimes \alpha_{0})_{\cal C}.
$$
We assign to ${\cal A}$ a bicomplex whose $n$-th column is 
$$
\prod_{[X_i]}\Bigl({\cal A}(X_0, X_1)[1] 
\otimes \ldots \otimes {\cal A}(X_{n-1}, X_n)[1]\otimes {\cal A}(X_{n}, X_0)[1]\Bigr)_{\cal C}.
$$
The differentials are induced by the differentials and the composition maps
 on Hom's. 
The cyclic homology complex ${\rm CC}_*({\cal A})$ of ${\cal A}$ 
is the total complex of this bicomplex. 
Its homology are the  cyclic homology 
of  ${\cal A}$. 

\vskip 2mm

Assume 
that  there are functorial pairings
$$
{\cal A}(X,Y)[1] \otimes {\cal A}(Y, X)[1] \lra {\cal H}^*.
$$ 
Then there is a morphism of complexes
\be \la{MofC}
{\rm HC}^*({\cal A})^* \lra {\rm CC}_*({\cal A}) \otimes {\cal H}.
\ee

For the category of harmonic bundles ${\rm Har}_X$ there is 
such a pairing with 
$$
{\cal H}:= H_{2n}(X)[-2].
$$
It provides a map
\be \la{6.27.11.4}
\varphi: {\rm Hom}\Bigl(H_0({\rm CC}_*({\rm Har}_X) \otimes 
{\cal H}, \C\Bigr)\lra {\rm HH}^0({\rm Har}_X)  \stackrel{(\ref{6.27.11.3})}{=}
H^0{\rm Fun}_{A_{\infty}}({\rm Har}_X, {\rm Har}_X).
\ee

\subsection{The Hodge correlators}

\bt \la{10.13.09.1}
a) There is a linear map, the 
Hodge correlator map
\be \la{HCM1}
{\rm Cor}_{\rm Har_X}: H_0({\rm CC}_*({\rm Har}_X) \otimes {\cal H}) \lra \C.
\ee
Combining it with  (\ref{6.27.11.4}), we get 
a cohomology class
\be \la{10.13.09.2}
{\bf  H}_{\rm Har_X} := \varphi({\rm Cor}_{\rm Har_X})\in 
H^0{\rm Fun}_{A_{\infty}}({\rm Har}_X, {\rm Har}_X).
\ee

b) There is a homotopy action of the 
 twistor Galois group ${G}_{\rm T}$ by $A_{\infty}$-autoequivalences of the category ${\rm Har}_X$ 
such that 

\begin{itemize}

\item
Its restriction to the subgroup $\C^*$ is the Simpson action  
\cite{Si1} 
on the category ${\rm Har}_X$. 

\item 
Its restriction to the Lie algebra ${\rm Lie}{\rm U}_{\rm Hod}$ 
is given by a Lie algebra map
\be \la{10.13.09.288}
{\Bbb H}_{\rm Har_X}:   {\rm L}_{\rm Hod}\lra 
{\rm H}^0{\rm Fun}_{A_{\infty}}({\rm Har}_X, {\rm Har}_X), 
\ee
uniquely determined by the condition that 
$
{\Bbb H}_{\rm Har_X}(G)  = {\bf  H}_{\rm Har_X}.
$ 
\end{itemize}

c) The action of the group ${G}_{\rm T}$ 
is functorial with respect to the pull backs.
\et

\subsection{Construction.} 

To define the Hodge correlator map
(\ref{HCM1}), we define a collection of degree zero maps
\be \la{HCM11}
{\rm Cor}_{\rm Hod_X}: 
\Bigl(H^\bullet(X, V^*_0\otimes V_1)[1] \otimes  \ldots \otimes H^\bullet(X, 
{V_m^* \otimes V_0})[1]\Bigr)_{\cal C} \otimes {\cal H} \lra \C.
\ee
The definition depends on some choices, like  
harmonic representatatives of cohomology classes. 
We prove that it is well defined on  $HC^0$, 
i.e. its resctriction to cycles is  independent of the choices, 
and coboundaries are mapped to zero.  

\begin{figure}[ht]
\centerline{\epsfbox{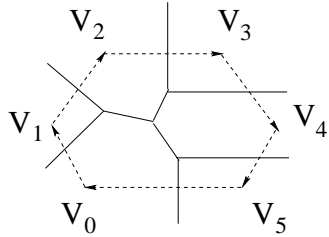}}
\caption{A decorated plane  trivalent tree; $V_i$ are harmonic bundles.}
\label{ht2}
\end{figure}

We picture an element in 
the sourse of the map (\ref{HCM11}) by a polygon $P$, see  Fig \ref{ht2}, 
whose vertices are 
the objects $V_i$, and the oriented  sides $V_iV_{i+1}$ 
are graded vector space
$
{\rm Ext}^*(V_i, V_{i+1})(1).
$

\paragraph{Green currents for harmonic bundles.} 
Let  $V$ be a harmonic bundle on $X$.  
Then there is a Doulbeaut  bicomplex 
$({\cal A}^{\bullet}(X, V); D', D'')$ where the differentials 
$D', D''$ are provided by the complex structure on $X$ and the harmonic metric on $V$. 
It satisfies the $D', D''$-lemma. 

Choose a splitting of the corresponding de Rham complex 
${\cal A}^\bullet(X, V)$ into 
an arbitrary subspace ${\cal H}ar^\bullet(X, V)$ isomorphically projecting onto 
the cohomology ${H}^\bullet(X, V)$ ("harmonic forms'') 
and its orthogonal complement. 
If $V = \C_X$, we choose  $a \in X$ and take the $\delta$-function $\delta_a$ 
at the point $a \in X$ 
as a representative of the fundamental class.

Let $\delta_{\Delta}$ be the Schwarz kernel of the identity map $V \to V$ given by 
the $\delta$-function of the diagonal, and 
$P_{\rm Har}$ the Schwarz kernel of the projector onto the space ${\cal H}ar^\bullet(X, V)$, 
realized by an $(n,n)$-form on $X \times X$. 
Choose a basis $\{\alpha_i\}$ 
in ${\cal H}ar^\bullet(X, V)$. Denote by 
$\{\alpha^{\vee}_i\}$ the dual basis. Then we have 
$$
P_{\rm Har} = \sum \alpha^{\vee}_i \otimes 
\alpha_i, \qquad \int_X  \alpha_i \wedge 
\alpha_j^{\vee} = \delta_{ij}.
$$
 Let $p_i:X \times X \to X$ be the projections onto the factors.  
\bd
A Green current $G(V; x,y)$ is a $p_1^*V^*\otimes p_2^*V$-valued current on $X \times X$, 
$$
G(V; x,y) \in {\cal D}^{2n-2}(X \times X, p_1^*V^*\otimes p_2^*V), \quad n={\rm dim}_\C X, 
$$
which satisfies the differential equation 
\be \la{GEQ}
(2\pi i)^{-1}D''D' G(V; x,y) = \delta_{\Delta} - P_{\rm Har}.
\ee
\ed
The two currents on the right hand side of (\ref{GEQ}) represent 
the same cohomology class, so the equation has a solution by the 
$D''D'$-lemma. 

\vskip 3mm
{\bf Remark}. The Green current depends on the choice of the ``harmonic forms''. 
So if $V = \C$, it depends on the choise of the base point $a$. 
Solutions of equation (\ref{GEQ}) are well defined modulo  
 ${\rm Im}D'' + {\rm Im}D' + {\cal H}ar^\bullet(X, V)$.

\paragraph{Construction of the Hodge correlators.} 
{\it Trees.} Take a plane trivalent tree 
$T$ dual to a triangulation of the polygon $P$, see Fig \ref{ht2}. 
The complement to $T$ in the polygon $P$ is a union of connected 
domains parametrized by the vertices of $P$, and thus 
decorated by the harmonic bundles $V_i$. Each edge $E$
of the tree $T$ is shared by two domains. 
The corresponding harmonic bundles are denoted 
$V_{E-}$ and $V_{E+}$. If $E$ is an external edge, we  
assume that $V_{E-}$ is before $V_{E+}$ 
for the clockwise orientation.

Given an internal vertex $v$ of the tree $T$, there are three domains 
sharing the vertex. We denote the corresponding harmonic bundles by $V_i, V_j, V_k$, 
where the cyclic order of the bundles 
agrees with the clockwise orientation. 
There is a natural trace map
\be \la{10.5.09.1}
{\rm Tr}_v: V_{i}^*\otimes V_{j} \otimes V_{j}^*\otimes V_{k} \otimes V_{k}^*\otimes V_{i}  = 
\lra \C.
\ee
It is invariant under the cyclic shift.

\vskip 3mm 
{\it Decorations}. For every edge $E$ of  $T$, 
choose a graded splitting of the de Rham complex 
$$
{\cal A}^\bullet(X, V_{E-}^*\otimes V_{E+}) = {\cal H}ar^\bullet(X, V_{E-}^*\otimes V_{E+}) \bigoplus 
{\cal H}ar^\bullet(X, V_{E-}^*\otimes V_{E+})^\perp.
$$ 
Then a decomposable class in
$\Bigl(\otimes_{i=0}^mH^*(X, V_i^\bullet\otimes V_{i+1})[1]\Bigr)_{\cal C}$ 
has a harmonic representative
$$
W= \Bigl(\alpha_{0,1} \otimes \alpha_{1, 2}\otimes \ldots 
\otimes {\alpha_{m, 0}}\Bigr)_{\cal C}.
$$
We are going to  assign to $W$ 
a  top degree current ${\kappa}(W)$ 
on 
\be \la{4.30.08.1}
X ^{\{\mbox{\rm internal vertices of $T$}\}}.
\ee
Each external edge $E$ of the tree $T$ is decorated by an element 
$$
\alpha_E \in {\cal H}ar^\bullet(X, V_{E-}^*\otimes V_{E+}).
$$
 Put the current $\alpha_E$ 
 to the copy of $X$ assigned  
to the internal vertex of the edge $E$, and pull it back to  (\ref{4.30.08.1}) using the  
projection $p_{\alpha_E}$ of the latter to the $X$. Abusing notation, we denote 
the pull back by $\alpha_E$. It is a form 
on (\ref{4.30.08.1}) with values in the bundle $p^*_{\alpha_E}(V_{E-}^*\otimes V_{E+})$

\begin{figure}[ht]
\centerline{\epsfbox{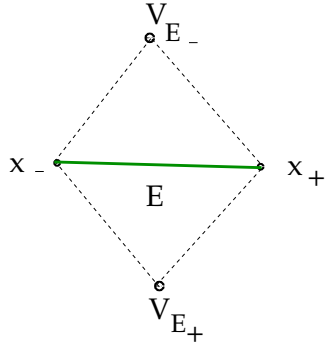}}
\caption{Decorations of the Green current assigned to an edge $E$.}
\label{ht3}
\end{figure}

\vskip 3mm
{\it Green currents}. 
We assign to each internal edge $E$ of the tree $T$ a Green current 
\be \la{10.5.09.3}
G(V_{E-}^*\otimes V_{E+}; x_-, x_+).
\ee
The order of $(x_-, x_+)$ agrees with the one of $(V_{E-}^*, V_{E+})$ as on
 Fig \ref{ht3}: the cyclic order  of  
$(V_{E-}^*, x_-, V_{E+}, x_+)$ 
agrees with the clockwise orientation. 
The Green current (\ref{10.5.09.3}) is symmetric: 
\be \la{10.5.09.4}
G(V_{E-}^*\otimes V_{E+}; x_-, x_+) = G(V_{E+}^*\otimes V_{E-}; x_+, x_-).
\ee
So it does not depend on the choice of orientation of the edge $E$.

 \vskip 3mm
{\it The map $\xi$}. 
There is a degree zero map 
\be \la{xixi}
\xi: {\cal A}^\bullet(X, V_0)[-1]\otimes \ldots  \otimes {\cal A}^\bullet(X, V_m)[-1] \lra
 {\cal A}^\bullet(X,V_0 \otimes ... \otimes V_m)[-1];
 \ee
\be \la{xixi1}
\varphi_0 \otimes \ldots \otimes
\varphi_m 
\lms {\rm Sym}_{\{0, ..., m\}}\Bigl(\varphi_0 \wedge D^\C \varphi_1 \wedge \ldots  \wedge D^\C 
\varphi_m\Bigr). 
\ee
The graded symmetrization in (\ref{xixi1}) is defined via isomorphisms 
$V_{\sigma(0)} \otimes ... \otimes V_{\sigma(m)}\to V_0 \otimes ... \otimes V_m $, 
where $\sigma$ is a permutation of $\{0, ..., m\}$.   
It is essential that ${\rm deg}D^\C \varphi = {\rm deg}\varphi +1$.

\vskip 3mm
{\it An outline of the construction}. 
We apply the operator $\xi$  to the product of 
the Green currents assigned to  
the internal edges of $T$. Then we  multiply  on (\ref{4.30.08.1})
the obtained local system valued current with the one provided by the decoration $W$,  
with an appropriate sign. Applying the product of the trace maps (\ref{10.5.09.1}) 
over the internal vertices of $T$, we get a top degree 
scalar current on (\ref{4.30.08.1}). 
Integrating it   we 
get a number assigned to $T$. Taking the sum over all plane trivalent 
trees $T$ decorated by $W$, we get 
a complex number ${\rm Cor}_{{\rm Har_X}}(W \otimes {\cal H})$. Altogether, we 
get the map (\ref{HCM1}). One checks that its degree is zero. 
The  signs in this 
definition are defined the same way as in \cite{G2}.

\bt
The  maps (\ref{HCM11}) give rise to a well defined  Hodge correlator map (\ref{HCM1}). 
\et

\vskip 3mm
\paragraph{Acknowledgments.} I am grateful to 
Alexander Beilinson and Maxim  Kontsevich for 
their interest to this project and useful discussions. 
This work was supported by the NSF grant DMS-1059129, MPI (Bonn) and IHES.

\end{document}